\documentclass{article}
\usepackage{amsmath}
\usepackage{amsfonts}
\usepackage{amssymb}

\newcommand{\norm}[1]{\left\lVert#1\right\rVert}
\newcommand{\abs}[1]{\left|#1\right|}

\newcommand{\Nz}{r}
\newcommand{\Nzw}{r'}
\renewcommand{\vec}[1]{\boldsymbol{#1}}
\newcommand{\lat}[1]{\mathcal{#1}}
\newtheorem{theorem}{Theorem}

\begin{document}

\title{On the approximation of real powers of sparse, infinite, bounded and Hermitian matrices}
\author{Roman Werpachowski\\Center for Theoretical Physics, Al. Lotnik\'{o}w 32/46\\02-668 Warszawa, Poland\\\texttt{roman@cft.edu.pl}}
\date{11 October 2006}

\maketitle

\begin{abstract}
We describe a method to approximate the matrix elements of a real power $\alpha$ of a positive (for $\alpha \in \mathbb{R}$) or non-negative (for $\alpha \ge 0$), infinite, bounded, sparse and Hermitian matrix $W$. The approximation uses only a finite part of the matrix $W$.

\hfill

\end{abstract}

\section{Introduction}

The motivation for this paper has been the research, carried out by Jerzy Kijowski~\cite{kijowski1,kijowski2,kijowski3} and, for the last year, the author, on the quantum field theory on a lattice. This physical theory arises from the discretization of the classical field Lagrangian and its subsequent quantization. As the number of degrees of freedom is infinite, Stone--von~Neumann's theorem~\cite{stone,vonneumann}, which proves the uniqueness (up to a unitary transformation) of the representation of the algebra of field operators, does not apply. The problem of finding such a representation has been customarily solved by resorting to global symmetries, such as the translational invariance. Our research has been concerned, however, with the quantum field theory on the background of non-homogeneous gravitational field. This eliminates global symmetries and renders traditional ways of finding the physical representation obsolete.

An attempt was made to find a representation of the field operators on a Hilbert space \textit{via} a limit of the representations on finite lattices, which are, as proven by von Neumann, uniquely defined up to a unitary transformation. When applied to a quadratic Hamiltonian on a $d$-dimensional rectangular lattice:
\[
\frac{1}{2} \sum_{\vec{n} \in \mathbb{Z}^d} \hat{\pi}_{\vec{n}}^2 + \sum_{\vec{m} \in \mathbb{Z}^d} \sum_{\vec{n} \in \mathbb{Z}^d} W_{\vec{m}\vec{n}} \hat{\phi}_{\vec{m}} \hat{\phi}_{\vec{n}} \ ,
\]
(where $\vec{n}$ is a vector of integer indices numbering the lattice sites, $\hat{\pi}_{\vec{n}}$ is the (real) momentum operator for the site $\vec{n}$, $\hat{\phi}_{\vec{n}}$ is the field operator, also for the site $\vec{n}$, and $W$ an infinite, positive, real and symmetric matrix), it turns out that the existence of the limit of representations is connected with the existence of the following limit: for each finite lattice $\lat{T} \subset \mathbb{Z}^d$, we define a finite matrix $W_{\lat{T}}$, which is a finite part of $W$ corresponding to the lattice $\lat{T}$. In order to make the problem well-posed physically, we need to impose some boundary conditions on the elements of $W_{\lat{T}}$ corresponding to the edges of $\lat{T}$. To ensure the existence of the groundstate, we choose non-negative boundary conditions. It turns out that the existence of the limit of finite-dimensional representations is guaranteed by the convergence of the sequence of matrix elements $(W^{\pm 1/2}_{\lat{T}})_{\vec{m}\vec{n}}$ with increasing lattices $\lat{T}$. Thus, the problem is a subset of a more general problem: for a real power of an infinite, positive and bounded Hermitian matrix $W$, we try to approximate its matrix elements with the corresponding elements of its finite parts, with possible boundary conditions added. This is the topic of this paper.

\section{Main result}

An infinite matrix $W$ is also a linear operator on the Hilbert space $\ell^2$ of square-summable infinite sequences. The matrix elements $W_{mn}$ are equal to the scalar product $(W e_n| e_m)$, where the $e_n$ are the canonical basis vectors of $\ell^2$. Conversely, any bounded linear operator on $\ell^2$ corresponds to a certain infinite matrix.

The infinite matrix $W^\alpha$ is, like $W$, a positive, bounded, linear and Hermitian operator on $\ell^2$ which, when acting on any eigenvector $v$ of $W$ with an eigenvalue $\lambda$, fulfills the equation $W^\alpha v = \lambda^\alpha v$. For a general theory of fractional powers of operators, see~\cite{carracedo}.

The core result of this paper is
\begin{theorem}
\label{theorem}
Consider an infinite Hermitian matrix $W$ with the following properties:
\begin{enumerate}
\item $\norm{W} < \infty$,
\item each row of $W$ has at most $k$ non-zero elements,
\item $\inf \sigma(W) \ge 0$, where $\sigma(W)$ is the spectrum of $W$ as an operator on $\ell^2$.
\end{enumerate}

For $\alpha \in \mathbb{R}$ (for $\inf \sigma(W) > 0$) or $\alpha \ge 0$ (for $\inf \sigma(W) = 0$), we have for a fixed pair of indices $m,n \in \mathbb{Z}$
\begin{equation}
\label{eq:limWPQalpha}
\lim_{P,Q \to \infty} (W_{PQ}^\alpha)_{mn} = W^\alpha_{mn} \quad \mathrm{(pointwise).}
\end{equation}
where $W_{PQ}$ is a $(P+Q+1) \times (P+Q+1)$ matrix with indices in the range $[-P,Q]$ which is defined as follows:
\begin{equation}
\label{eq:defWPQint}
(W_{PQ})_{mn} = W_{mn} \ , \quad -P < m < Q \vee -P < n < Q
\end{equation}
and
\begin{equation}
\label{eq:defWPQbnd}
(W_{PQ})_{mn} = W_{mn} + D_{mn} \ , \quad m,n \in \lbrace -P, Q \rbrace \ .
\end{equation}
The Hermitian matrix $D$ fulfills the following conditions: for any $P,Q > 0$,
% and for any $x,y \in \mathbb{C}$, 
we have
\begin{equation}
\label{eq:Dcond0}
W_{PQ} \ge \inf \sigma(W)
\end{equation}
(i.e. the lowest eigenvalue of $W_{PQ}$ is greater or equal to $\inf \sigma(W)$) and
\begin{equation}
\label{eq:Dcond}
%0 \le D_{-P,-P} \abs{x}^2 + D_{QQ} \abs{y}^2 + D_{-P,Q} \cc{x} y + D_{Q,-P} x \cc{y} < d \ , \quad 0 \le d < \infty \ .
\norm{W_{PQ}} \le \norm{W} + d \ ,
\end{equation}
for $d \ge 0$.
\end{theorem}

\section{Physical example}

Consider an infinite matrix $M$ which generates the following quadratic form for $x \in \ell^2$:
\[
x^\dagger M x = a \sum_{n=-\infty}^\infty \abs{x_n}^2 + b \sum_{n=-\infty}^\infty \abs{x_n - x_{n-1}}^2 \ , \quad a > 0 \, \wedge \, b > 0 \ .
\]
It arises from the homogeneous theory of scalar field on a one-dimensional lattice. The finite matrix $M_{PQ}$ is chosen in such a way that it generates the following quadratic form for $x\in \mathbb{C}^{P+Q+1}$:
\[
x^\dagger M_{PQ} x = a \sum_{x=-P}^Q \abs{x_n}^2 + b \sum_{n=1-P}^Q \abs{x_n - x_{n-1}}^2 + b\abs{x_Q - x_{-P}}^2 \ ,
\]
an example of the so-called ``periodic boundary conditions'' in physics. It turns out that for $-P \le m,n \le Q$ and $\alpha \in \mathbb{R}$ we have
\[
(M_{PQ}^\alpha)_{mn} = \sum_{k=0}^{P+Q} \frac{\left[ a + b \sin^2 \left( \frac{\pi k}{P+Q+1} \right) \right]^\alpha}{P+Q+1} \cos \left( \frac{2 \pi k (m - n)}{P+Q+1} \right) \ .
\]
On the other hand,
\[
M^\alpha_{mn} = \int_0^1 (a + b \sin^2 \pi k)^\alpha \cos [2 \pi k (m - n)] \mathrm{d} k \ .
\]
It is clear that
\[
\lim_{P,Q} (M_{PQ}^\alpha)_{mn} = M^\alpha_{mn} \ .
\]

\section{Proof of Theorem~\ref{theorem}}

For $\alpha = 0$, the theorem is trivial.

Since every row of $W$ has at most $k$ non-zero elements, any matrix element of any positive and integer power of $W$ can be calculated as a finite sum of finite products of the matrix elements of $W$:
\begin{equation}
\label{eq:Wpowers}
\begin{split}
W^2_{mn} &= \sum_{m_1 \in \Nz(m)} W_{mm_1} W_{m_1n} \\
W^3_{mn} &= \sum_{m_1 \in \Nz(m)} \sum_{m_2 \in \Nz(m_1)} W_{mm_1} W_{m_1m_2} W_{m_2n} \\
W^4_{mn} &= \dots \ ,
\end{split} 
\end{equation}
where
\[
\Nz(m) := \lbrace n | W_{mn} \neq 0 \rbrace \ .
\]
Because of~\eqref{eq:defWPQint} and~\eqref{eq:Wpowers}, for any indices $m,n$ and for every non-negative and integer power $j$ one can choose such finite $P_j$ and $Q_j$ that 
\[
\forall_{P\ge P_j,Q\ge Q_j} \  W^j_{mn} = (W_{PQ}^j)_{mn} \ .
\]
This proves~\eqref{eq:limWPQalpha} for $\alpha \in \mathbb{N}$.

Because of~\eqref{eq:defWPQint},~\eqref{eq:defWPQbnd},~\eqref{eq:Dcond0} and~\eqref{eq:Dcond}, we have the following bounds for the spectrum of $W_{PQ}$ \emph{as an operator on $\mathbb{C}_{PQ}$, the space of all complex finite sequences indexed from $-P$ to $Q$}:
\[
%\label{eq:WPQbnd}
\sigma_{PQ}(W_{PQ}) \subset [ \inf \sigma(W), \norm{W} + d ] \ .
\]
We make the distinction between $\sigma_{PQ}$ (spectrum of an operator on $\mathbb{C}_{PQ}$) and $\sigma$ (spectrum of an operator on $\ell^2$), because $W_{PQ}$, as an operator on $\ell^2$, has infinitely many eigenvectors with the zero eigenvalue:
\[
W_{PQ} e_{n} = 0 \ , \quad n < -P \, \vee \, n > Q \ ,
\]
which means that
\[
\sigma(W_{PQ}) = \sigma_{PQ}(W_{PQ}) \cup \{ 0 \} \ .
\]

The lower and upper bound will be denoted by $c$ and $w$, respectively:
\[
%\label{eq:cw}
\begin{split}
c &:= \inf \sigma(W) \ , \quad c \ge 0 \\
w &:= \norm{W} + d \ , \quad c \le w < \infty \ .
\end{split}
\]
Since $\norm{W} \le w$, $W^\alpha$ can be expressed~\cite{maurin} as the following power series:
\begin{equation}
\label{eq:Walpha}
W^\alpha = w^\alpha \sum_{j=0}^\infty {\alpha \choose j} \left( \frac{W - w I}{w} \right)^j \ ,
\end{equation}
where $I$ is the infinite identity matrix, $I_{mn} = \delta_{m,n}$. We know that $\norm{W - w I} \le w - c \le w$. Hence, $\norm{W - wI}/w \le 1$ and the series~\eqref{eq:Walpha} is weakly convergent (in the case of $\norm{W - wI} / w = 1$ the series is weakly convergent for $\alpha \ge 0$ only). Weak convergence of~\eqref{eq:Walpha} implies the convergence of the following series:
\begin{equation}
\label{eq:Walphamn}
(W^\alpha)_{mn} = w^\alpha \sum_{j=0}^\infty {\alpha \choose j} \frac{(W - w I)^j_{mn}}{w^j} \ ,
\end{equation}
where $(W^\alpha)_{mn}$ is defined as a matrix element of the linear operator $W^\alpha$ between two canonical basis vector of $\ell^2$, $e_m$ and $e_n$.

Since every row of $W$ has at most $k$ non-zero elements, any matrix element of any positive and integer power of $W - w I$ can be calculated as a finite sum of finite products of the matrix elements of $W - w I$:
\begin{equation}
\label{eq:powers}
\begin{split}
(W - w I)^2_{mn} &= \sum_{m_1 \in \Nzw(m)} (W - w I)_{mm_1} (W - w I)_{m_1n} \\
(W - w I)^3_{mn} &= \sum_{m_1 \in \Nzw(m)} \sum_{m_2 \in \Nzw(m_1)} (W - w I)_{mm_1} (W - w I)_{m_1m_2} (W - w I)_{m_2n} \\
(W - w I)^4_{mn} &= \dots \ ,
\end{split} 
\end{equation}
where
\[
\Nzw(m) := \lbrace n | (W - w I)_{mn} \neq 0 \rbrace \ .
\]
Because of~\eqref{eq:defWPQint} and~\eqref{eq:powers}, we can define a non-negative integer $j_{PQ}$ such that
\begin{equation}
\label{eq:jPQdef}
j < j_{PQ} \ \Leftrightarrow \  (W - w I)^j_{mn} = (W_{PQ} - w I_{PQ})^j_{mn} \ ,
\end{equation}
where $I_{PQ}$ is the $(P+Q+1) \times (P+Q+1)$ identity matrix, and
\begin{equation}
\label{eq:limjPQ}
\lim_{P,Q \to \infty} j_{PQ} = \infty \ .
\end{equation}

$(W_{PQ}^\alpha)_{mn}$, for $W_{PQ}$ being an operator on $\mathbb{C}_{PQ}$, can be expressed as a series similar to~\eqref{eq:Walphamn}:
\[
%\label{eq:WPQalphamn}
(W^\alpha_{PQ})_{mn} = w^\alpha \sum_{j=0}^\infty {\alpha \choose j} \frac{(W_{PQ} - w I_{PQ})^j_{mn}}{w^j} \ .
\]
We split this series with the help of~\eqref{eq:jPQdef}:
\begin{equation}
\label{eq:WPQalphapmn}
(W^\alpha_{PQ})_{mn} = w^\alpha \sum_{j=0}^{j_{PQ} - 1} {\alpha \choose j} \frac{(W - w I)^j_{mn}}{w^j} + w^\alpha \sum_{j=j_{PQ}}^\infty {\alpha \choose j} \frac{(W_{PQ} - w I_{PQ})^j_{mn}}{w^j} \ .
\end{equation}
Using~\eqref{eq:Walphamn} and~\eqref{eq:WPQalphapmn}, we get
\begin{multline*}
(W^\alpha_{PQ})_{mn} - (W^\alpha)_{mn} =\\= w^\alpha \left[ \sum_{j=j_{PQ}}^\infty {\alpha \choose j} \frac{(W_{PQ} - w I_{PQ})^j_{mn}}{w^j} - \sum_{j=j_{PQ}}^\infty {\alpha \choose j} \frac{(W - w I)^j_{mn}}{w^j} \right] \ .
\end{multline*}
The difference can be estimated as follows:
\begin{multline}
\label{eq:diffbnd}
\abs{(W^\alpha_{PQ})_{mn} - (W^\alpha)_{mn}} = \\
w^\alpha \abs{\sum_{j=j_{PQ}}^\infty {\alpha \choose j} \frac{(W_{PQ} - w I_{PQ})^j_{mn}}{w^j} - \sum_{j=j_{PQ}}^\infty {\alpha \choose j} \frac{(W - w I)^j_{mn}}{w^j}} \\
\le w^\alpha \left[ \sum_{j=j_{PQ}}^\infty \abs{{\alpha \choose j}} \frac{\abs{(W_{PQ} - w I_{PQ})^j_{mn}}}{w^j} + \sum_{j=j_{PQ}}^\infty \abs{{\alpha \choose j}} \frac{\abs{(W - w I)^j_{mn}}}{w^j} \right] \ .
\end{multline}
It is a well-known fact that, for every Hermitian operator $A$ on a Hilbert space $\mathcal{H}$ which has a countable orthonormal basis $\{ e_n \}$, every matrix element of $A$ fulfills the inequality
\[
\abs{A_{mn}} \le \norm{A} \ .
\]
Indeed, using the Cauchy-Schwarz inequality we obtain
\[
\abs{A_{mn}} = \abs{( A e_n | e_m  )} \le \norm{e_m} \norm{A e_n} \le \norm{A} \norm{e_m} \norm{e_n} = \norm{A} \ .
\]
Since $\ell^2$ has a countable orthonormal basis and $W - w I$ is Hermitian, we have
\[
\abs{(W - w I)^j_{mn}} \le (w - c)^j \ .
\]
Ditto
\[
\abs{(W_{PQ} - w I_{PQ})^j_{mn}} \le (w - c)^j \ ,
\]
since $W_{PQ} - w I_{PQ}$ is a Hermitian matrix. We insert these bounds into~\eqref{eq:diffbnd} and obtain
\begin{equation}
\label{eq:diffbnd2}
\abs{(W^\alpha_{PQ})_{mn} - (W^\alpha)_{mn}} < 2 w^\alpha \sum_{j=j_{PQ}}^\infty \abs{{\alpha \choose j}} \left( \frac{w - c}{w} \right)^j \ .
\end{equation}

The series 
\begin{equation}
\label{eq:diffser}
\sum_{j=0}^\infty \abs{{\alpha \choose j}} \left( \frac{w - c}{w} \right)^j
\end{equation}
(i.e. the same as above with $j_{PQ}$ set to 0) converges (as with~\eqref{eq:Walpha}, for $c = 0$ it converges only for $\alpha \ge 0$).
% For non-negative and integer $\alpha$, we have
% \begin{equation}
% \label{eq:sum1}
% \sum_{j=0}^\infty \abs{{\alpha \choose j}} \left( \frac{w - c}{w} \right)^j = \sum_{j=0}^{\infty} {\alpha \choose j} \left( \frac{w - c}{w} \right)^j = (2 - c/w)^\alpha \ ,\quad \alpha \in \mathbb{Z}^+ \ .
% \end{equation}
For negative $\alpha$ and $c > 0$,
\begin{equation}
\label{eq:sum2}
\sum_{j=0}^\infty \abs{{\alpha \choose j}} \left( \frac{w - c}{w} \right)^j = \sum_{j=0}^{\infty} {\alpha \choose j} \left( \frac{c - w}{w} \right)^j = (c/w)^\alpha\ , \quad \alpha < 0 \ .
\end{equation}
For positive and non-integer $\alpha$, we have ($[\alpha]$ denotes the integer part of $\alpha$):
\begin{multline}
\label{eq:sum3}
\sum_{j=0}^\infty \abs{{\alpha \choose j}} \left( \frac{w - c}{w} \right)^j = \\
% \sum_{j=0}^{[\alpha]+1} {\alpha \choose j} \left( \frac{w - c}{w} \right)^j +\\
% + (-1)^{[\alpha]+1} \sum_{j=[\alpha]+2}^\infty {\alpha \choose j} \left( \frac{c - w}{w} \right)^j = \sum_{j=0}^{[\alpha]+1} {\alpha \choose j} \left( \frac{w - c}{w} \right)^j -\\
% - (-1)^{[\alpha]+1} \sum_{j=0}^{[\alpha]+1} {\alpha \choose j} \left( \frac{c - w}{w} \right)^j + (-1)^{[\alpha]+1} \sum_{j=0}^\infty {\alpha \choose j} \left( \frac{c - w}{w} \right)^j = \\
% = \sum_{j=0}^{[\alpha]+1} {\alpha \choose j} \frac{(w-c)^j + (-1)^{[\alpha]} (c-w)^j}{w^j}  +  (-1)^{[\alpha]+1} \sum_{j=0}^\infty {\alpha \choose j} \left( \frac{c - w}{w} \right)^j =\\
= \sum_{j=0}^{[\alpha]+1} {\alpha \choose j} \frac{(w-c)^j + (-1)^{[\alpha]} (c-w)^j}{w^j} + (-1)^{[\alpha]+1} (c/w)^\alpha \ , \quad \alpha > 0 \wedge \alpha \notin \mathbb{Z} \ .
\end{multline} 
Convergence of~\eqref{eq:diffser}, proven in~\eqref{eq:sum2} and~\eqref{eq:sum3}, shows that the sum of the elements of the series~\eqref{eq:diffbnd} from $j = j_{PQ}$ to $\infty$ goes to zero for $j_{PQ} \to \infty$:
\begin{equation}
\label{eq:limtail}
\lim_{j_{PQ} \to \infty} 2 w^\alpha \sum_{j=j_{PQ}}^\infty \abs{{\alpha \choose j}} \left( \frac{w - c}{w} \right)^j = 0 \ .
\end{equation}

Using~\eqref{eq:limjPQ},~\eqref{eq:diffbnd2} and~\eqref{eq:limtail}, we calculate:
\begin{multline*}
\lim_{P,Q \to \infty} \abs{(W^\alpha_{PQ})_{mn} - (W^\alpha)_{mn}} < \lim_{P,Q \to \infty} 2 w^\alpha \sum_{j=j_{PQ}}^\infty \abs{{\alpha \choose j}} \left( \frac{w - c}{w} \right)^j \\
= \lim_{j_{PQ} \to \infty}  2 w^\alpha \sum_{j=j_{PQ}}^\infty \abs{{\alpha \choose j}} \left( \frac{w - c}{w} \right)^j = 0 \ ,
\end{multline*}
% Since $\abs{(W^\alpha_{PQ})_{mn} - (W^\alpha)_{mn}} \ge 0$, we finally obtain
% \[
% \lim_{P,Q \to \infty} \abs{(W^\alpha_{PQ})_{mn} - (W^\alpha)_{mn}} = 0 \ ,
% \]
which proves~\eqref{eq:limWPQalpha}.

\section{Accuracy and applicability of the approximation}

Inequality~\eqref{eq:diffbnd2} and results~\eqref{eq:sum2} and~\eqref{eq:sum3} can be used to estimate the accuracy of our approximation. For a given unknown matrix element $W^\alpha_{mn}$, we are interested in how fast the right side of~\eqref{eq:diffbnd2} diminishes with increasing $P$ and $Q$. This involves two problems:
\begin{enumerate}
\item the relation between $P,Q$ and $j_{PQ}$, which is determined by the structure of $W$ near the chosen matrix element, and
\item the speed of convergence of the series~\eqref{eq:diffser}, which appears in~\eqref{eq:diffbnd2}.
\end{enumerate}
Since the series~\eqref{eq:diffser} has only positive elements, its convergence is determined by how fast they vanish with increasing $j$. This in turn is determined by the power factor $(1-c/w)^j$ (for $c > 0$) and the binomial coefficient $\alpha \choose j$. The smaller is $1 - c/w$, the faster will~\eqref{eq:diffser} converge. The dependence of $\alpha \choose j$ on $j$ is such that~\eqref{eq:diffser} will converge more slowly when $\alpha \lesssim -1$, than in the opposite case.

The relation between $P,Q$ and $j_{PQ}$ (i.e. how fast $j_{PQ}$ diverges when $P,Q \to \infty$) cannot be determined without additional assumptions about $W$. It is obvious that $j_{PQ} = 1$ (the minimal value it can take) for $-P > \min(m,n)$ or $Q < \max(m,n)$, where $m,n$ are indices of the matrix element that we want to compute. 

Suppose, for example, that $W$ is $k$-diagonal for $k = 2l+1$, so that $W_{mn} = 0$ for $\abs{m - n} > l$ and is non-zero otherwise. We then have $j_{PQ} = 1 + [\min(\min(m,n) + P - 1, Q - \max(m,n) - 1)/l]$ for $-P \le \min(m,n) - 1$ and $Q \ge \max(m,n) + 1$. Asymptotically, $j_{PQ}$ grows linearly with $\min(P,Q)$.

Since $W$ is sparse, positive integer powers may be easily computed without resorting to any approximations, using~\eqref{eq:Wpowers}.

Our approximation can be used to approximate locally the solutions to infinite matrix equations $W x = f$ (i.e. approximate the solutions of linear operator equations, see~\cite{bakusinskij,chen3}), for sequences $f$ which are non-zero only for a finite subset of indices. Additionally, the approximation might be applied also to calculate numerically the matrix element of some power of a finite, but very large matrix which fulfills the conditions of our theorem.

\section{Summary and generalization}

We have shown that one can use a finite part of the infinite, non-negative, Hermitian, sparse bounded matrix $W$ to approximate a chosen matrix element of a power $W^\alpha$ for any real or positive $\alpha$ (for $\inf W > 0$ or $\inf W = 0$, respectively). The accuracy of the approximation depends on the bounds of the matrix spectrum, the value of $\alpha$ and on the structure of $W$ near the chosen element. Possible applications of the approximation have been discussed briefly.

It is apparent that the result presented in this paper can be easily generalized to the case of not only power functions, but also all functions which are finite and possess a power series expansion absolutely convergent on the whole spectrum of $W$.

\section{Acknowledgements}

The author wishes to thank prof. Jerzy Kijowski for inspiring this work and dr Andrzej Wakulicz for encouragment to the writing of this paper. The author also acknowledges the support of the Polish Ministry of Science and Higher Education under the grant ``Quantum Information and Quantum Engineering'' No.~PBZ-MIN-008/P03/2003.

\bibliography{apprpow}

\providecommand{\bysame}{\leavevmode\hbox to3em{\hrulefill}\thinspace}
\providecommand{\MR}{\relax\ifhmode\unskip\space\fi MR }
\providecommand{\MRhref}[2]{%
  \href{http://www.ams.org/mathscinet-getitem?mr=#1}{#2}
}
\providecommand{\href}[2]{#2}
\begin{thebibliography}{1}

\bibitem{bakusinskij}
Anatolij~Borisovi\v{c} Baku\v{s}inskij and M.~Yu. Kokurin, \emph{Iterative
  methods for approximate solution of inverse problems}, Springer, 2004.

\bibitem{carracedo}
Celso~Mart\'{i}nez Carracedo and Miguel~Sanz Alix, \emph{The theory of
  fractional powers of operators}, Elsevier, 2001.

\bibitem{chen3}
Mingjun Chen, Zhongying Chen, and Guanrong Chen, \emph{Approximate solutions of
  operator equations}, World Scientific, 1997.

\bibitem{kijowski3}
Jerzy Kijowski, Gerd Rudolph, and Cezary \'{S}liwa, \emph{Charge superselection
  sectors for scalar {QED} on the lattice}, Ann. Henri Poincar\'{e} \textbf{4}
  (2003), no.~6, 1137--1167.

\bibitem{kijowski2}
Jerzy Kijowski, Gerd Rudolph, and Artur Thielmann, \emph{Algebra of observables
  and charge superselection sectors for {QED} on the lattice}, Comm. Math.
  Phys. \textbf{188} (1997), no.~3, 535--564.

\bibitem{kijowski1}
Jerzy Kijowski and Artur Thielmann, \emph{Quantum electrodynamics on a
  space-time lattice}, J. Geom. Phys. \textbf{19} (1996), 173--205.

\bibitem{maurin}
Krzysztof Maurin, \emph{Methods of {Hilbert} spaces}, PWN, Warsaw, 1967.

\bibitem{stone}
M.~H. Stone, \emph{Linear transformations in {H}ilbert space, {III}:
  {O}perational methods and group theory}, Proc. Nat. Acad. Sci. U.S.A.
  \textbf{16} (1930), 172--175.

\bibitem{vonneumann}
J.~von Neumann, \emph{Die {E}indeutigkeit der {S}chr\"{o}dingerschen
  {O}peratoren}, Math. Ann. \textbf{104} (1931), 570--578.

\end{thebibliography}
\bibliographystyle{amsplain}

\end{document}